\documentclass[conference,a4paper]{IEEEtran}
\IEEEoverridecommandlockouts
\usepackage{graphics}
\usepackage[demo]{graphicx}
\usepackage{epsfig,amssymb}
\usepackage{amsmath,mathtools}
\usepackage{times}
\usepackage{mathrsfs}
\usepackage{array}
\usepackage{amsmath, latexsym, amsfonts, amssymb}
\usepackage[mathscr]{eucal}
\usepackage{array, tabularx}
\usepackage[table]{xcolor}
\usepackage{psfrag}
\usepackage{multirow}\usepackage{booktabs}
\usepackage{epsfig}
\usepackage{cite}
\usepackage{tikz}
\usepackage{scalerel}
\usepackage{url}
\usepackage{xcolor}
\usepackage[english]{babel}
\usepackage[autostyle]{csquotes}
\usepackage{authblk}
\usepackage{subcaption}
\usepackage{mwe}
\usepackage{algorithm}
 \usepackage[]{algpseudocode}

\newcommand{\brparen}[1]{\left\{#1\right\}}

\usepackage{graphics}
\usepackage[demo]{graphicx}
\usepackage{epsfig,amssymb}
\usepackage{amsmath,mathtools}
\usepackage{times}
\usepackage{mathrsfs}
\usepackage{array}
\usepackage{amsmath, latexsym, amsfonts, amssymb}
\usepackage[mathscr]{eucal}
\usepackage{array, tabularx}
\usepackage[table]{xcolor}
\usepackage{psfrag}
\usepackage{multirow}\usepackage{booktabs}
\usepackage{epsfig}
\usepackage{cite}
\usepackage{tikz}
\usepackage{scalerel}
\usepackage{url}
\usepackage{xcolor}
\usepackage[english]{babel}
\usepackage[autostyle]{csquotes}
\usepackage{authblk}
\usepackage{subcaption}
\usepackage{adjustbox}

\newcommand{\ie}{\ensuremath{{\text{\em i.e.}}}}

\def\BibTeX{{\rm B\kern-.05em{\sc i\kern-.025em b}\kern-.08em
    T\kern-.1667em\lower.7ex\hbox{E}\kern-.125emX}}
\begin{document}

\title{Performance Enhancement of C-V2X Mode 4 Utilizing Multiple Candidate Single-subframe Resources\hspace{-10pt}}


\author{Geeth~P.~Wijesiri N.B.A,~\IEEEmembership{Member,~IEEE,}
        Tharaka~Samarasinghe,~\IEEEmembership{Senior Member,~IEEE,}
        and~Jussi Haapola,\\~\IEEEmembership{Member,~IEEE,}
\thanks{G. P. Wijesiri N.B.A is with the Department of Electronic and Telecommunication Engineering, University of Moratuwa, Sri Lanka, and the Department of Electrical and Information Engineering, University of Ruhuna, Sri Lanka (e-mail: geeth@eie.ruh.ac.lk).

T. Samarasinghe is with the Department of Electronic and Telecommunication Engineering, University of Moratuwa, Sri Lanka, and the Department of Electronic and Electrical  Engineering, University of Melbourne,
Australia (e-mail: tharakas@uom.lk).

J. Haapola is with the Centre for Wireless Communications, University of Oulu, Finland (e-mail: jussi.haapola@oulu.fi).
}
\thanks{This research has partially received funding from the European Union’s Horizon 2020 research and innovation programme under grant agreement No 857031 (project 5G!Drones), the Academy of Finland 6Genesis Flagship (grant 318927) and the AHEAD/RA3/RIC/MRT/ITS-Multidisciplinary Transport Development Project of the University of Moratuwa under World Bank grant 6026-LK/8743-LK.}}

\bibliographystyle{ieeetr}
\maketitle


\begin{abstract}
Prioritization of data streams in cellular vehicle-to-everything (C-V2X) may lead to unfavorable packet delays in low priority streams. This paper studies the allocation of multiple candidate single-subframe resources (CSRs) per vehicle as a solution. It proposes a methodology to determine the number of CSRs for each vehicle based on the number of total vehicles, and to assign the multiple data streams among them for simultaneous transmission. The numerical results highlight the achievable delay gains of the proposed approach, and its negligible impact on packet collisions.
\end{abstract}

\section{Introduction}
In Release 14, the third generation partnership project (3GPP) introduced cellular vehicle–to–everything (C-V2X) Mode 4 to support vehicular communications without the support of cellular infrastructure. Therein, the medium access control (MAC) layer plays a crucial role in handling stringent, but varying, delay and reliability constraints of different V2X applications. The variable delay constraints have necessitated the technologies to support multi-priority data streams \cite{b17}. 
Recently in \cite{jp2}, it has been shown that in the presence of multi-priority data streams, the competing technology IEEE 802.11p \cite{ieee_802_11_2016} outperforms C-V2X Mode 4 in terms of delay and priority management, thanks to its enhanced distributed channel access (EDCA) mechanism. 
Thus, addressing this issue is of importance for C-V2X Mode 4, specifically in avoiding stale packets in lower priority data streams.

Limited work on enhancing the performance of C-V2X Mode 4 can be found in the literature. Performance enhancement through variations in the transmit power is studied in \cite{c_v2x_en_3}. The work in \cite{c_v2x_en_1} and \cite{omnetpp_sim} are more related, and they focus on spectrum management and the semi-persistence scheduling (SPS) algorithm. In C-V2X Mode 4, vehicles use the SPS algorithm in a distributed manner to sense the radio resources (called candidate single-subframe resources (CSRs) utilized by other vehicles in a sensing window, and to select a CSR for its own transmission. 
To this end, \cite{c_v2x_en_1} introduces a weighted power averaging methodology for sensing the CSRs in the sensing window. The authors of \cite{omnetpp_sim} show that re-selecting the same CSR used for the previous transmission (reusing) more frequently, and using exponential sensing window sizes under high channel load levels can lead to enhanced performance. 
Both works limit their study to a single data stream at each vehicle. Moreover, as the performance metric, they focus mainly on the packet reception ratio of C-V2X Mode 4, which is known to be superior compared to its competing technology in \cite{jp2}. The authors of \cite{jp2} further highlight the lower channel utilization exhibited by C-V2X Mode 4, which we try to exploit to ameliorate the pressing concerns on delay and priority management.

In this paper, we focus on multi-priority data streams, and study the possibility of enhancing the performance of C-V2X Mode 4 by allocating multiple CSRs for each vehicle, which in turn increases the channel utilization. 
The procedure leads to two fundamental problems.
\begin{itemize}
\item \textbf{Determining the number of CSRs per vehicle: } \textcolor{black}{We determine how many CSRs can be allocated to each vehicle based on the total number of vehicles in the system, without coordination.}
\item \textbf{Allocating the multi-priority data streams among the CSRs: } Having established the number of CSRs, we introduce a procedure for allocating the multi-priority data streams among the allocated CSRs. 
\end{itemize}

We evaluate the performance of the proposed method using the discrete-time Markov chain (DTMC)-based models in \cite{jp2}.
The results show that the parallelism achieved by allocating multiple CSRs leads to significant reductions in the average delay, specially in the low priority data streams. In general, this is achieved by allocating separate CSRs for low priority data streams and allowing them more frequent transmission opportunities, opposed to waiting till all higher priority queues are empty. We can observe that the optimal group selection depends on the number of CSRs available and the generation rates of each data stream. It should be noted that increased packet collisions is the tradeoff of allocating multiple CSRs. However, since C-V2X Mode 4 inherently has superior collision resolution as compared to IEEE 802.11p \cite{jp2}, the increase is incremental, and can be considered insignificant when compared with the perceived benefits in terms of delay and priority management.

The paper is organized as follows. The system model is presented in Section \ref{sec:secII}. Section \ref{sec:secIII} studies the allocation of multiple CSRs to a vehicle. Section \ref{sec:secIV} presents the numerical results and the discussion, and Section \ref{sec:secV} concludes the paper.

\section{System Model}\label{sec:secII}
We consider a network with $N$ vehicles, and each vehicle transmits decentralized environmental notification messages (DENM), high priority DENM (HPD), multi-hop DENM (MHD), and cooperative awareness messages (CAM).  The priority order of serving these packets according to the standard is as follows: HPD $>$ DENM $>$ CAM $>$ MHD. We use subscripts $i \in \mathcal{I} = \brparen{H,D,C,M}$ to differentiate between the parameters for HPD, DENM, CAM, and MHD, respectively. The CAM packets are generated periodically with a generation interval $T_{C}$. 
HPD, DENM, and MHD are randomly generated at an average generation rate of $\lambda_{m}$, for $m\in \mathcal{I} \setminus\{C\}$, based on events initiated by human activity or environmental conditions. Thus, the parameters related to packet generation can be written as $\mathcal{P} =\brparen{\lambda_{H},\lambda_{D},\lambda_{M},T_{C}}$.
The generated packets are queued separately, and transmitted according to their level of priority, based on the SPS algorithm. The network allocates each vehicle $n_{CSR}$ CSRs for this transmission, and the network allocates the multi-priority data streams among the $n_{CSR}$ CSRs for possible enhanced performance at the MAC layer. 
To this end, CSRs are adjacent sub-channel sets within the subframe that are large enough to fit in the sidelink control information (SCI) and the transport block (TB) to be transmitted.

\begin{figure}[t]
\begin{center}

\tikzset{every picture/.style={line width=0.75pt}} 

\begin{tikzpicture}[x=0.75pt,y=0.75pt,yscale=-1,xscale=1]

\draw   (299,104.29) .. controls (299,95.12) and (306.43,87.69) .. (315.6,87.69) -- (383.4,87.69) .. controls (392.57,87.69) and (400,95.12) .. (400,104.29) -- (400,154.09) .. controls (400,163.26) and (392.57,170.69) .. (383.4,170.69) -- (315.6,170.69) .. controls (306.43,170.69) and (299,163.26) .. (299,154.09) -- cycle ;
\draw    (250,160) -- (296,159.7) ;
\draw [shift={(298,159.69)}, rotate = 539.63] [color={rgb, 255:red, 0; green, 0; blue, 0 }  ][line width=0.75]    (10.93,-3.29) .. controls (6.95,-1.4) and (3.31,-0.3) .. (0,0) .. controls (3.31,0.3) and (6.95,1.4) .. (10.93,3.29)   ;
\draw    (398,99) -- (444,98.7) ;
\draw [shift={(446,98.69)}, rotate = 539.63] [color={rgb, 255:red, 0; green, 0; blue, 0 }  ][line width=0.75]    (10.93,-3.29) .. controls (6.95,-1.4) and (3.31,-0.3) .. (0,0) .. controls (3.31,0.3) and (6.95,1.4) .. (10.93,3.29)   ;
\draw    (399,130) -- (445,129.7) ;
\draw [shift={(447,129.69)}, rotate = 539.63] [color={rgb, 255:red, 0; green, 0; blue, 0 }  ][line width=0.75]    (10.93,-3.29) .. controls (6.95,-1.4) and (3.31,-0.3) .. (0,0) .. controls (3.31,0.3) and (6.95,1.4) .. (10.93,3.29)   ;
\draw    (398,161.69) -- (444,161.39) ;
\draw [shift={(446,161.38)}, rotate = 539.63] [color={rgb, 255:red, 0; green, 0; blue, 0 }  ][line width=0.75]    (10.93,-3.29) .. controls (6.95,-1.4) and (3.31,-0.3) .. (0,0) .. controls (3.31,0.3) and (6.95,1.4) .. (10.93,3.29)   ;
\draw    (253,99) -- (299,98.7) ;
\draw [shift={(301,98.69)}, rotate = 539.63] [color={rgb, 255:red, 0; green, 0; blue, 0 }  ][line width=0.75]    (10.93,-3.29) .. controls (6.95,-1.4) and (3.31,-0.3) .. (0,0) .. controls (3.31,0.3) and (6.95,1.4) .. (10.93,3.29)   ;
\draw    (252,130) -- (298,129.7) ;
\draw [shift={(300,129.69)}, rotate = 539.63] [color={rgb, 255:red, 0; green, 0; blue, 0 }  ][line width=0.75]    (10.93,-3.29) .. controls (6.95,-1.4) and (3.31,-0.3) .. (0,0) .. controls (3.31,0.3) and (6.95,1.4) .. (10.93,3.29)   ;

\draw (264,146) node [anchor=north west][inner sep=0.75pt]    {$A$};
\draw (402,80) node [anchor=north west][inner sep=0.75pt]    {$d_{a}{}_{v}{}_{g,l}( A) ,\ l\in B$};
\draw (403,112) node [anchor=north west][inner sep=0.75pt]    {$P_{c}{}_{o}{}_{l}( A)$};
\draw (404,145) node [anchor=north west][inner sep=0.75pt]    {$CU( A)$};
\draw (306,116) node [anchor=north west][inner sep=0.75pt]   [align=left] {\begin{minipage}[lt]{63.39pt}\setlength\topsep{0pt}
\begin{center}
DTMC based\\overall model
\end{center}
 
\end{minipage}};
\draw (267,85) node [anchor=north west][inner sep=0.75pt]    {$N$};
\draw (266,116) node [anchor=north west][inner sep=0.75pt]    {$\Gamma $};

\end{tikzpicture}
\end{center}
\caption{The DTMC based overall model for a given CSR based on the modeling techniques in \cite{jp2}.}
\label{fig:overall_model}
\end{figure}
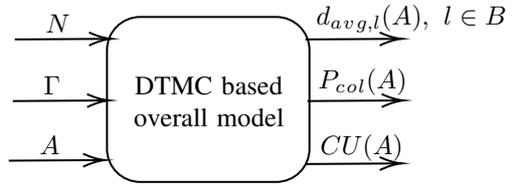

In \cite{jp2}, a similar setup is modeled using DTMCs for $n_{CSR}=1$. The overall model in \cite{jp2} consists of four separate DTMCs that model the generation of the packet types of interest, four DTMCs that model their device-level packet queues, and one DTMC that models the MAC layer operations related to the transmission. We directly resort to the modeling techniques and the derived performance measures in \cite{jp2} to evaluate the performance achieved when 
$n_{CSR}\geq 1$.
For a given CSR, $A \subset \mathcal{P}$, $N$, and $\Gamma$, which is called the selection window size in the SPS algorithm, act as the inputs to the model, as illustrated in Fig. \ref{fig:overall_model}. $A$ is based on the data streams allocated to that particular CSR, and the respective performance metrics which are functions of $A$ act as the outputs of the model. To this end, we obtain the average delay of the $l$-th data stream $d_{avg,l}(A)$, for each $l\in B$, the collision probability $P_{col}(A)$, and the channel utilization $CU(A)$ as outputs, where $B \subset \mathcal{I}$ is the corresponding set of indices for $A$, {\em e.g.}, $A= \brparen{\lambda_H,T_c} \ \rightarrow \ B= \brparen{H,C}$. 

\begin{figure}[t]
    \centering
    \includegraphics[scale=.53]{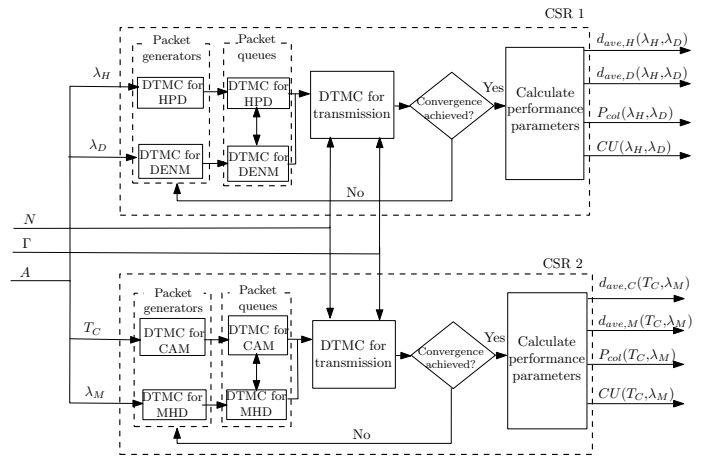}
    \caption{A diagram illustrating the overall model for an example scenario where $n_{CSR}=2$.}
    \label{fig:overrall_model_eg}
 \end{figure}

We use an example scenario to further elaborate the inputs, outputs, and the usage of the DTMCs. Consider an example where two CSRs are used at each vehicle. Moreover, HPD and DENM data are allocated for the first CSR, and the other two data streams are allocated for the second, as illustrated in Fig. \ref{fig:overrall_model_eg}. With regards to the first CSR, $N$, $\Gamma$, and $A=\brparen{\lambda_H,\lambda_D}$ act as the inputs, and $d_{avg,H}(\lambda_H,\lambda_D)$, $d_{avg,D}(\lambda_H,\lambda_D)$, $P_{col}(\lambda_H,\lambda_D)$ and $CU(\lambda_H,\lambda_D)$ act as the outputs. Further on the DTMCs, two DTMCs model the generation of HPD and DENM packets, and two DTMCs model their device-level packet queues. The priority management is incorporated in the queue models such that a resultant queue is connected with the next DTMC that models the transmission of these packets. The dependencies among these DTMCs are appropriately modelled according to \cite{jp2}, and illustrated in Fig. \ref{fig:overrall_model_eg}. We obtain expressions for the steady state probabilities of the DTMCs and solve them iteratively until convergence is achieved. Upon convergence, the performance metrics are calculated using the steady state probabilities, which act as the outputs. A similar procedure is followed for the second CSR.    
It is not hard to see that the network performance depends on the number of CSRs used by each vehicle, and how the data streams are allocated among the CSRs. This is the main focus of the next section.

\section{Allocation of Multiple CSRs to a Vehicle}\label{sec:secIII}

Allocating multiple CSRs to a vehicle leads to two key fundamental questions. Firstly, we need to ascertain  a plausible value for $n_{CSR}$. We assume $n_{CSR} \leq 4$ for simplicity, $\ie$, we assume that multiple CSRs are not allocated to a single data stream.  Secondly, given $n_{CSR}$, we need to decide how the four parallel data streams should be allocated among the $n_{CSR}$ CSRs. In this section, we find solutions to these two questions based on the average delay, and present an algorithm for multiple CSR allocation for a vehicle.

\subsection{Determining the Number of CSRs}

It is rather intuitive that $n_{CSR}$ is inversely proportional to $N$. Also, $N$ has a direct impact on the length of the selection window $\Gamma$, which is defined as the maximum latency in ms \cite{b4}, and should be set at either 20 ms, 50 ms or 100 ms according to the standard \cite{b2}. Increasing $\Gamma$ increases the time gap between two successive transmissions, which in turn leads to more radio resources for the network. Thus, $CSR_{tot}$ is proportional to $\Gamma$, where $CSR_{tot}$ denotes the total number of available CSRs in the selection window. The standard allows allocating 80\% of these CSRs to the users. Thus, $N_{\max}= 0.8 CSR_{tot} / n_{CSR}$, where $N_{\max}$ is the maximum number of users that can be handled simultaneously. With larger $\Gamma$ values, the system can support more users, but with a trade-off of a higher delay. For a given $N$, we need to obtain $n_{CSR}$, and then appropriately set $\Gamma$ values for each of these CSRs. 

To this end, we resort to allocating the same $\Gamma$ value across the multiple CSRs allocated to a vehicle. This allocation can be justified as follows. Although different $\Gamma$ values are allocated across the multiple CSRs, the number of simultaneous users the overall system can support is constrained by the shortest $\Gamma$ value currently being used in the system. For example, consider a scenario where each vehicle uses 2 CSRs, and we allocate $\Gamma=$20 ms ($N_{\max}=400$) for the first CSR and $\Gamma=$50 ms ($N_{\max}=1000$) for the second. It is not hard to see that the overall system can only support 400 users. Thus, allocating different $\Gamma$ values across multiple CSRs is counter productive, and only leads to higher delay values at the MAC layer. 

\begin{figure}[t]
    \centering
    \includegraphics[scale=.55]{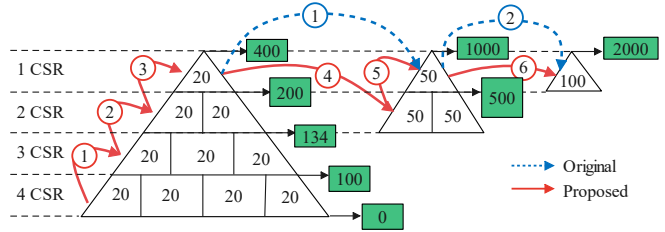}
    \caption{The process of selecting $n_{CSR}$ and $\Gamma$ with $N$.}
    \label{fig:csr_sw}
 \end{figure}

The relationship between $n_{CSR}$ and $N_{\max}$ for the three values of $\Gamma$ is tabulated in Table \ref{tab:maximum_neighbours_vs_CSRs}. We prefer shorter $\Gamma$ values to reduce the average delay \cite{jp2}, and higher values for $n_{CSR}$ to exploit higher degrees of freedom (resources) for the allocation of the parallel data streams. Fig. \ref{fig:csr_sw} illustrates the selection of $n_{CSR}$ and $\Gamma$, when $N$ increases from 1 to 2000. The methodology associated with allocating multiple CSRs is presented using the solid red arrows, and the dashed blue arrows show the equivalent transitions when the vehicle uses a single CSR as per the standard. The values inside the rectangles of the pyramidal shapes depict the value of $\Gamma$, and the values inside the green colored rectangles depict the maximum supported value of $N$, for each multiple CSR configuration.

\begin{table}[ht]
\caption{The relationship between $n_{CSR}$ and $N_{\max}$ for different values of $\Gamma$. }
\label{tab:maximum_neighbours_vs_CSRs}
\centering
\begin{tabular}{|c|c|c|c|}
\hline
\multirow{2}{*}{$n_{CSR}$} & \multicolumn{3}{c|}{\textbf{\begin{tabular}[c]{@{}l@{}}$N_{\max}$ \end{tabular}}} \\ \cline{2-4} 
                              & \textbf{$\Gamma$=20 ms}                           & \textbf{$\Gamma$=50 ms}                           & \textbf{$\Gamma$=100 ms}                          \\ \hline
1                             & 400                                       & 1000                                     &2000                                     \\ \hline
2                             & 200                                      & 500                                       & 1000                                     \\ \hline
3                             & 134                                      & 334                                      & 667                                      \\ \hline
4                             & 100                                     & 250                                      & 500                                      \\ \hline
\end{tabular}
\end{table}

\subsection{Allocating the Multi-priority Data Streams Among the CSRs}

Having established $n_{CSR}$ and an appropriate value for $\Gamma$ for a given $N$, we now face the problem of allocating the multi-priority data streams among the CSRs. We present the grouping options using Table \ref{tab:grouping_option}. Let $n_G$ denote the index of the grouping option. If $n_{CSR}=1$, we simply allocate all data streams to the single CSR, as per the current standard, and we call this grouping option $n_{G}=1$. If there are 4 CSRs, the grouping is again trivial, as we can allocate a separate CSR for each data stream. Thus, we get $n_G=15$. 
The number of groups for the 2 and 3 CSR scenarios is based on the Stirling number of the second kind, thus there are 7 and 6 grouping options for each scenario, respectively.

\begin{table}[t]
\caption{Grouping options with associated generation parameters.}
 \label{tab:grouping_option}
\centering
\begin{tabular}{|c|c|c|c|c|c|}
\hline
\textbf{$n_{csr}$}              & \textbf{$n_{G}$} & \multicolumn{4}{c|}{\textbf{Generation Parameters for CSRs}}                                                          \\ \hline
\multirow{2}{*}{\textbf{1 CSR}} & \textbf{}        & \multicolumn{4}{c|}{\textbf{CSR1}}                                                                                    \\ \cline{2-6} 
                                & 1                & \multicolumn{4}{c|}{$\lambda_{H}$, $\lambda_{D}$, $\lambda_{M}$, $T_{C}$}                                             \\ \hline
\multirow{8}{*}{\textbf{2 CSRs}} & \textbf{} & \multicolumn{2}{c|}{\textbf{CSR 1}} & \multicolumn{2}{c|}{\textbf{CSR 2}} \\ \cline{2-6} 
                                & 2                & \multicolumn{2}{c|}{$T_{C}$}                      & \multicolumn{2}{c|}{$\lambda_{H}$,  $\lambda_{D}$, $\lambda_{M}$} \\ \cline{2-6} 
                                & 3                & \multicolumn{2}{c|}{$\lambda_{H}$}                & \multicolumn{2}{c|}{$\lambda_{D}$, $T_{C}$, $\lambda_{M}$}        \\ \cline{2-6} 
                                & 4                & \multicolumn{2}{c|}{$\lambda_{D}$}                & \multicolumn{2}{c|}{$\lambda_{H}$, $T_{C}$, $\lambda_{M}$}        \\ \cline{2-6} 
                                & 5                & \multicolumn{2}{c|}{$\lambda_{M}$}                & \multicolumn{2}{c|}{$\lambda_{H}$, $T_{C}$, $\lambda_{D}$}        \\ \cline{2-6} 
                                & 6                & \multicolumn{2}{c|}{$\lambda_{H}$, $T_{C}$}       & \multicolumn{2}{c|}{$\lambda_{D}$, $\lambda_{M}$}                 \\ \cline{2-6} 
                                & 7                & \multicolumn{2}{c|}{$\lambda_{H}$, $\lambda_{M}$} & \multicolumn{2}{c|}{$\lambda_{D}$, $T_{C}$}                      \\ \cline{2-6} 
                                & 8                & \multicolumn{2}{c|}{$\lambda_{H}$, $\lambda_{D}$} & \multicolumn{2}{c|}{$T_{C}$, $\lambda_{M}$}                       \\ \hline
\multirow{7}{*}{\textbf{3 CSRs}}         &                  & \textbf{CSR 1}           & \textbf{CSR 2}         & \multicolumn{2}{c|}{\textbf{CSR 3}}                               \\ \cline{2-6} 
                                & 9                & $T_{C}$                  & $\lambda_{H}$          & \multicolumn{2}{c|}{$\lambda_{D}$, $\lambda_{M}$}                 \\ \cline{2-6} 
                                & 10               & $\lambda_{H}$            & $\lambda_{D}$          & \multicolumn{2}{c|}{$T_{C}$, $\lambda_{M}$}                       \\ \cline{2-6} 
                                & 11               & $T_{C}$                  & $\lambda_{M}$          & \multicolumn{2}{c|}{$\lambda_{H}$, $\lambda_{D}$}                 \\ \cline{2-6} 
                                & 12               & $\lambda_{D}$            & $T_{C}$               & \multicolumn{2}{c|}{$\lambda_{H}$, $\lambda_{M}$}                  \\ \cline{2-6} 
                                & 13               & $\lambda_{D}$            & $\lambda_{M}$          & \multicolumn{2}{c|}{$\lambda_{H}$, $T_{C}$}                        \\ \cline{2-6} 
                                & 14               & $\lambda_{H}$            & $\lambda_{M}$          & \multicolumn{2}{c|}{$\lambda_{D}$, $T_{C}$}                       \\ \hline
\multirow{2}{*}{\textbf{4 CSRs}} & \textbf{} & \textbf{CSR 1}   & \textbf{CSR 2}   & \textbf{CSR 3}   & \textbf{CSR 4}   \\ \cline{2-6} 
                                & 15               & $\lambda _{H}$           & $\lambda_{D}$          & $T_{C}$                       & $\lambda_{M}$                     \\ \hline
\end{tabular}
\end{table}

We select the best grouping option $n_{G}^{\star}$ for a given $n_{CSR}$ with respect to the average delay. Let $\mathcal{D}_{n_{G},l}$ be the average delay of the $l$-th data stream for grouping option $n_{G}$. These delay values can be calculated using the DTMC models by appropriately setting the parameter combinations tabulated in Table \ref{tab:grouping_option} for $A$, as explained in reference to Fig. \ref{fig:overrall_model_eg}.  For example, the resulting average delay values for grouping option $n_{G}=2$ are $\mathcal{D}_{2,C}=d_{avg,C}(T_{C})$ and $\mathcal{D}_{2,m}=d_{avg,m}(\lambda_{H}, \lambda_{D}, \lambda_{M})$ for $m \in \mathcal{I}\setminus \{C\}$. The sum average delay for grouping option $n_{G}$ is written as
$$ \Delta \mathcal{D}_{n_{G}}=\sum_{l\in \mathcal{I}}  w_{l}\mathcal{D}_{n_{G},l},$$
where  $w_{l}$, for $l\in \mathcal{I}$, denotes a weight for each data stream of interest. We set $w_H > w_D > w_C > w_M$ such that the priority management in the standard \cite{b17} is incorporated in our selection, and we select the sum average delay minimizing grouping option as the best one for a given $n_{CSR}$.
Our approach of allocating multiple CSRs to a vehicle is formally presented through Algorithm \ref{alg:enhanced_v2x}. 

\begin{algorithm}[t]\caption{\small Multiple CSR Allocation for C-V2X Mode 4.}\label{alg:enhanced_v2x}
\begin{algorithmic}[1]
\small
\Procedure{$n_{CSR}$ \& $\Gamma$ allocation }{$N$,$\lambda_{H}$,$\lambda_{D}$,$\lambda_{M}$,$T_{C}$}
\State \hspace{-10pt}$\Gamma=0$, $n_{CSR}=0$ \Comment{Initialization}
\State \hspace{8pt}\textbf{if}\hspace{8pt} {$0<N\leq100$} \hspace{5.5pt}$\rightarrow$ $n_{CSR}=4$, $\Gamma=20$ ms 
\State \hspace{-10pt}\textbf{else if} {$100< N\leq134$} \hspace{5pt}$\rightarrow$
$n_{CSR}=3$, $\Gamma=20$ ms 
\State \hspace{-10pt}\textbf{else if} {$134<N\leq200$} \hspace{5pt}$\rightarrow$
 $n_{CSR}=2$, $\Gamma=20$ ms 
\State \hspace{-10pt}\textbf{else if} {$200<N\leq400$} \hspace{5pt}$\rightarrow$ 
$n_{CSR}=1$, $\Gamma=20$ ms
\State \hspace{-10pt}\textbf{else if} {$400<N\leq500$} \hspace{5pt}$\rightarrow$ 
$n_{CSR}=2$, $\Gamma=50$ ms 
\State \hspace{-10pt}\textbf{else if} {$500< N\leq1000$} $\rightarrow$ 
$n_{CSR}=1$, $\Gamma=50$ ms 
\State \hspace{-10pt}\textbf{end if}
\State \hspace{-10pt} BestGroup($n_{CSR}$,$\Gamma$,$\lambda_{H}$,$\lambda_{D}$,$\lambda_{M}$,$T_{C}$)
\EndProcedure

 \vspace{3pt}
  \State\hspace{20pt}\textbf{function} BestGroup($n_{CSR}$,$\Gamma$,$\lambda_{H}$,$\lambda_{D}$,$\lambda_{M}$,$T_{C}$)
  \State\hspace{25pt}$best\_group=0$ \Comment{Initialization}
     \State\hspace{43pt}\textbf{if} $(n_{CSR}=1)$ \hspace{-2pt}$\rightarrow$ $n_{G}^{\star}=1$
      \State\hspace{25pt}\textbf{else if} {$(n_{CSR}=2)$}\hspace{-1.5pt}
       $\rightarrow$ $n_{G}^{\star} = \underset{k\in\{2,\dots,8\}}{\arg\min}\  \Delta \mathcal{D}_{k}$
      \State\hspace{25pt}\textbf{else if} {$(n_{CSR}=3)$}\hspace{-1.5pt}
       $\rightarrow$ $n_{G}^{\star} = \underset{k\in\{9,\dots,14\}}{\arg\min}\  \Delta \mathcal{D}_{k}$
       \State\hspace{25pt}\textbf{else if} {($n_{CSR}=4)$} \hspace{-1pt}$\rightarrow$ $n_{G}^{\star}=15$

       \State\hspace{25pt}\textbf{end if} 
         \State\hspace{25pt}\textbf{Return} $n_{G}^{\star}$
  \State\hspace{20pt}\textbf{end function}
\end{algorithmic}
\end{algorithm}

\section{Numerical Results and Discussion}\label{sec:secIV}
This section presents numerical results that highlight the performance of using multiple CSRs by the vehicles according to Algorithm \ref{alg:enhanced_v2x}. 
The reference packet formats of HPD, DENM, CAM, and MHD are set according to \cite{CAM_standard, DENM_standard}.
HPD and DENM packets are retransmitted at fixed intervals for added reliability \cite{denm_rep} as per the standard, and for the results, we have set the number of repetitions to 8 and 5 times, respectively. 
The candidate values for $\lambda_{m}$, where $m \in \mathcal{I} \setminus\{C\}$, are set by being consistent with the example use case scenarios in \cite{etsi_qos} and  $T_{C}$ is set between 100 ms and 1 s according to standard \cite{CAM_standard}. We set $w_{H}=0.4$, $w_{D}=0.3$, $w_{C}=0.2$, and $w_{M}=0.1$ to extend the priority management in the standard \cite{b17} into the grouping methodology.



\subsection{Average Delay}

\begin{table*}[t]
\caption{The average delay reduction percentages with $N$ for CAM and MHD packets.}
\label{tab:average_delay_reduction}
\centering
\begin{tabular}{|c|c|c|c|c|c|c|c|c|c|c|c|c|c|c|}
\hline
\cellcolor[HTML]{FFFFFF}{\color[HTML]{333333} } &
   &
   &
  \multicolumn{12}{c|}{\textbf{The average delay reduction (\%)}} \\ \cline{4-15} 
\cellcolor[HTML]{FFFFFF}{\color[HTML]{333333} } &
   &
   &
  \multicolumn{6}{c|}{\textbf{$\lambda_{H}=\lambda_{D}=0.1$, $\lambda_{M}=1$ packets/s}} &
  \multicolumn{6}{c|}{\textbf{$\lambda_{H}=\lambda_{D}=1$, $\lambda_{M}=10$ packets/s}} \\ \cline{4-15} 
\cellcolor[HTML]{FFFFFF}{\color[HTML]{333333} } &
   &
   &
  \multicolumn{3}{c|}{\textbf{$T_{C}=100$ ms}} &
  \multicolumn{3}{c|}{\textbf{$T_{C}=500$ ms}} &
  \multicolumn{3}{c|}{\textbf{$T_{C}=100$ ms}} &
  \multicolumn{3}{c|}{\textbf{$T_{C}=500$ ms}} \\ \cline{4-15} 
\multirow{-4}{*}{\cellcolor[HTML]{FFFFFF}{\color[HTML]{333333} \textbf{$N$}}} &
  \multirow{-4}{*}{\textbf{$n_{CSR}$}} &
  \multirow{-4}{*}{\textbf{\begin{tabular}[c]{@{}c@{}}$\Gamma$\\ (ms)\end{tabular}}} &
  \textbf{$n_{G}^{\star}$} &
  \textbf{CAM} &
  \textbf{MHD} &
  \textbf{$n_{G}^{\star}$} &
  \textbf{CAM} &
  \textbf{MHD} &
  \textbf{$n_{G}^{\star}$} &
  \textbf{CAM} &
  \textbf{MHD} &
  \textbf{$n_{G}^{\star}$} &
  \textbf{CAM} &
  \textbf{MHD} \\ \hline
\textbf{$[0.100)$} &
  4 &
  20 &
  15 &
  9.8 &
  14.2 &
  15 &
  0.1 &
  0.6 &
  15 &
  11.5 &
  26.2 &
  15 &
  0.2 &
  16.3 \\ \hline
\textbf{$[100,134)$} &
  3 &
  20 &
  11 &
  9.8 &
  14.2 &
  14 &
  0.04 &
  0.6 &
  11 &
  11.5 &
  26,2 &
  14 &
  0.13 &
  16.3 \\ \hline
\textbf{$[134,200)$} &
  2 &
  20 &
  2 &
  9.8 &
  11.0 &
  5 &
  0.01 &
  0.6 &
  2 &
  11.5 &
  18.1 &
  5 &
  0.05 &
  16.3 \\ \hline
\textbf{$[200,400)$} &
  1 &
  20 &
  1 &
  - &
  - &
  1 &
  - &
  - &
  1 &
  - &
  - &
  1 &
  - &
  - \\ \hline
\textbf{$[400,500)$} &
  2 &
  50 &
  2 &
  54.7 &
  50.4 &
  5 &
  0.03 &
  3.2 &
  2 &
  77.4 &
  80.7 &
  5 &
  0.3 &
  73.1 \\ \hline
\end{tabular}
\end{table*}

\begin{figure}[t]
    \centering
    \includegraphics[scale=.37]{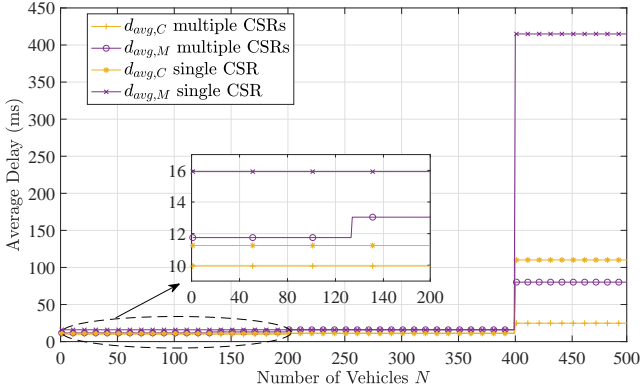}
    \caption{The behavior of average delay vs $N$ at $\lambda_H=\lambda_D=1$,  $\lambda_M=10$ packets/s and $T_{C}=100$ ms.}
    \label{fig:delay}
 \end{figure}

Firstly, we present the results on the average delay. Fig \ref{fig:delay} illustrates the average delay behavior of CAM and MHD packets with $N$ for $\lambda_H=\lambda_D=1$,  $\lambda_M=10$ packets/s and $T_{C}=100$ ms. 
As shown in this figure, the average delay is mainly sensitive to $\Gamma$, exhibiting a step-wise increase when $\Gamma$ switches from a shorter to a longer value, and constant with respect to $N$ for fixed $\Gamma$.
It was noticed that delay gains for HPD and DENM were negligibly small due to their higher priority levels, thus omitted in the results. The higher priority packets are served first regardless of the number of CSRs, and hence, the gains are insignificant. It can be seen that the average delay values can be maintained below 100 ms thanks to the utilization of multiple CSRs. 
Further results on the average delay reduction percentages for CAM and MHD data streams relative to using a single CSR as per the standard, are tabulated in Table  \ref{tab:average_delay_reduction}.

Firstly, while focusing on less critical scenarios such as roadwork warnings and safety function out of normal condition warnings, where the packet generation rates are considerably lower ($\lambda_D=0.1$ packets/s), we can observe clear gains of utilizing multiple CSRs at each vehicle. The maximum average delay reduction percentage for CAM is 54.7\%, which is around 30 ms, and for MHD, it is 50.4\%, which is around 48.9 ms. We can expect the gains for MHD to increase further at higher $\lambda_M$ values. For example, if $\lambda_M=10$, the gain is 69.2\% which is 236 ms.     
While focusing on more critical scenarios, such as emergency electronic brake lights and warnings from emergency vehicles that have higher packet generation rates ($\lambda_D=1$), we can observe very high gains for both CAM and MHD,  \textit{i.e.,} a maximum delay reduction percentage of 77.4\% (85 ms) and  80.7\% (334 ms), respectively, which can also be observed in the range $N \in [400,500]$ in Fig. \ref{fig:delay}. Thus, the multiple CSR configurations can contribute considerably to alleviate the issue of stale packets in low priority queues. In general, the results show that the proposed method works favorably for both less critical and critical scenarios when considering the average delay.
Furthermore, it is useful when more frequent location updates are required, which is achieved using CAM. This can be seen by comparing the gains for $T_{C}=100$ ms with $T_{C}=500$ ms. Thus, the multiple CSR configuration will be ideal for applications that require high CAM rates, which is in the magnitude of 10 Hz according to the standard \cite{etsi_qos}.

Moreover, based on the $n_{G}^{\star}$ values, we can provide the following insights on allocating the multiple CSRs among the data streams. For 2 CSR configurations, it can be observed that allocating the periodic and event-triggered traffic for separate CSRs performs better when the CAM rate is high. On the other hand, better delay-wise performance can be obtained by allocating a separate CSR for MHD when the CAM rate is low. This eliminates the necessity of MHD packets waiting till all higher priority queues are empty. For 3 CSR configurations, better delay-wise performance can be obtained by allocating a CSR each for CAM and MHD streams, and the other CSR for HPD and DENM streams, given that the CAM rate is high. As the system considers both HPD and DENM packets to have relatively higher priority, and hence transmits with minimum delay, allocating a CSR each for these two data streams does not show an advantage in this scenario. However, at lower CAM rates, we allocate a CSR each for MHD, HPD and DENM streams, and CAM can be grouped with the data stream having the lowest generation rate. In our results, we observe that the grouping was with DENM as it has a lower effective rate compared to HPD due to the lower number of packet repetitions. 

\subsection{Collision Probability and Channel Utilization}
\begin{figure}[t]
    \centering
    \includegraphics[scale=.37]{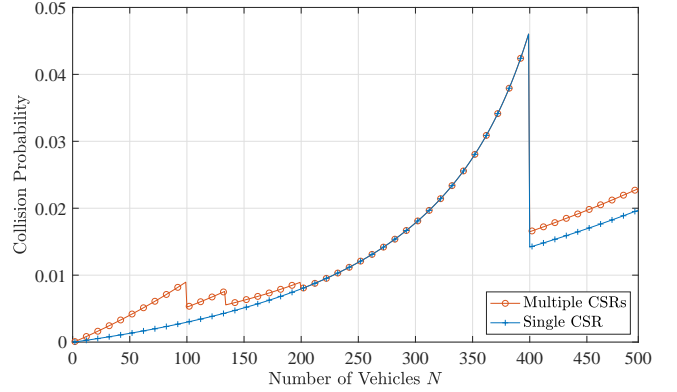}
    \caption{The behavior of collision probability vs $N$ at $\lambda_H=\lambda_D=0.1$, $\lambda_M=1$ packets/s and $T_{C}=100$ ms.}
    \label{fig:P_col}
 \end{figure}

Although favorable in terms of delay, using multiple CSRs may lead to a trade-off in terms of higher packet collisions. The variation of per CSR collision probability for the single and multiple CSR configurations are illustrated in Fig. \ref{fig:P_col}. 
The collision probability increases exponentially with $N$. We can also observe that $\Gamma$ has a significant impact on the collision probability as the collision probability reduces significantly when the value of $\Gamma$ increases. This is thanks to the availability of more radio resources at higher $\Gamma$ values \cite{jp2}. 

When comparing the collision probabilities of the two configurations, it can be seen that utilizing multiple CSRs lead to higher collision probability values in $0 < N\leq 200$  and $400 < N \leq 500$. This is due to the higher number of overlaps in the selection windows when multiple CSRs are used as explained in \cite{tcom_paper}. 
Therefore, there is a clear tradeoff of using multiple CSRs. However, it can be seen from  Fig. \ref{fig:P_col} that the maximum increase in collision probability is approximately 0.6\% (at $N=99$) when compared to using a single CSR, which is rather insignificant compared to the gains achieved on delay and priority management.  

\begin{figure}[t]
    \centering
    \includegraphics[scale=.37]{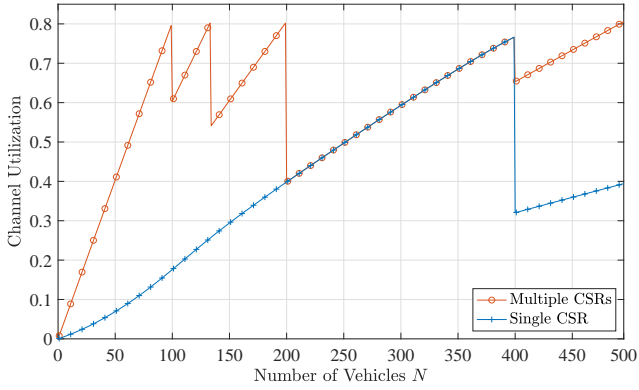}
    \caption{The behavior of channel utilization vs $N$ at $\lambda_H=\lambda_D=0.1$, $\lambda_M=1$ packets/s and $T_{C}=100$ ms.}
    \label{fig:cu}
 \end{figure}

We already saw that the average delay increases with $\Gamma$, and in Section \ref{sec:secIII} we showed that each $\Gamma$ has its respective $N_{max}$. We can increase the $N_{max}$ threshold levels further if the standard allows allocating a higher percentage of available CSRs to the users, $\ie$, increasing the 80\% parameter in the SPS algorithm stated in Section \ref{sec:secIII}. In that case, allocating multiple CSRs may deem even more favorable is terms of average delay. However, we can observe from Fig. \ref{fig:P_col} that this change leads to an exponential increase in the collision probability. Therefore, fine-tuning $N_{max}$ needs to be done only after carefully studying the QoS requirements of the applications. The behavior of channel utilization with $N$ is shown in Fig. \ref{fig:cu}. The figure clearly shows how the channel utilization has been exploited in the ranges of $0 < N\leq 200$  and $400 < N \leq 500$ to achieve the initial objectives.

We end the discussion by providing some insights on some implications of the proposed method. 
Firstly, let us focus on the SPS algorithm. In the SPS algorithm according the current standard, each vehicle is capable of tracking the CSRs used by itself and the neighboring vehicles within the sensing window. These identified CSRs are excluded when selecting a CSR for the subsequent transmission. Using multiple CSRs at each vehicle will not bring about major changes to how the SPS algorithm tracks CSRs used by neighboring vehicles. However, the SPS algorithm needs to be slightly modified to identify and exclude the CSRs used within the target vehicle itself, to minimize the internal collisions. The authors note that the proposed method may also cause changes in the hardware setup as parallel transmission is required, but with the developments in multi-antenna technologies, handling the hardware implications seems practically feasible. Extensive details on hardware implications are beyond the scope of this paper.

\section{Conclusions}\label{sec:secV}
This paper has focused on a vehicular network that utilizes C-V2X Mode 4 for communication, and supports multi-priority data streams to fulfill varying quality-of-service constraints of ITS use cases. It has proposed increasing the channel utilization of the network through the allocation of multiple CSRs at each vehicle, and has studied its achievable performance gains at the MAC layer. The proposed method has led to two fundamental questions, which are, how many CSRs should be allocated to each vehicle, and how the multi-priority data streams should be allocated among them. The number of CSRs at each vehicle has been ascertained as a function of total number of vehicles in the network, and a procedure has been introduced for allocating the multi-priority data streams among them based on the average delay in the network. The results have shown that using multiple CSRs at each vehicle can lead to significant gains in the network in terms of average delay. In particular, the average delay of lower priority data streams can be improved significantly by allocating them separate CSRs, which ameliorates the risk of stale packets. As a trade-off, an increase in the collision probability can be observed, but the performance loss is almost insignificant compared to the delay gains.

\ifCLASSOPTIONcaptionsoff
  \newpage
\fi

\bibliography{Geeth-bibfile}
\end{document}